\documentstyle[12pt]{article}

\input amssym.def
\input amssym.tex

\newcommand{\maps}{\mapsto}

\newcommand{\eop}{\bigstar}  

\newcommand{\card}[1]{{\vert #1 \vert} }

\newcommand{\otp}[1]{\hbox{otp($#1$)}}
\newcommand{\forces}{\Vdash}
\newcommand{\decides}{\parallel}

\newcommand{\dom}{{\rm dom}}
\newcommand{\Dom}{{\rm Dom}}

\newcommand{\rge}{{\rm rge}}
\newcommand{\crit}{{\rm crit}}

\newcommand{\cf}{{\rm cf}}

\newcommand{\implies}{\Longrightarrow}

\newtheorem{THEOREM}{Theorem}[section]

\newtheorem{Conclusion}[THEOREM]{Conclusion}

\newtheorem{LEMMA}[THEOREM]{Lemma}

\newtheorem{Main Theorem}[THEOREM]{Main Theorem}
\newenvironment{main Theorem}{\begin{Main Theorem}} 
{\end{Main Theorem}}

\newtheorem{Theorem}[THEOREM]{Theorem}

\newtheorem{Definition}[THEOREM]{Definition}

\newtheorem{Conventions}[THEOREM]{Conventions}

\newtheorem{Main Definition}[THEOREM]{Main Definition}
\newenvironment{main definition}{\begin{Main Definition}}
{\end{Main Definition}}

\newtheorem{Lemma}[THEOREM]{Lemma}

\newtheorem{Notation}[THEOREM]{Notation}

\newtheorem{Convention}[THEOREM]{Convention}

\newtheorem{Note}[THEOREM]{Note}

\newtheorem{Observation}[THEOREM]{Observation}

\newtheorem{Remark}[THEOREM]{Remark}

\newtheorem{Main Fact}[THEOREM]{Main Fact}
\newenvironment{main Fact}{\begin{Main Fact}}{\end{Main Fact}}

\newtheorem{Fact}[THEOREM]{Fact}

\newtheorem{Subfact}[THEOREM]{Subfact}

\newtheorem{Claim}[THEOREM]{Claim}

\newtheorem{Main Claim}[THEOREM]{Main Claim}
\newenvironment{main claim}{\begin{Main Claim}}{\end{Main Claim}}

\newtheorem{Corrolary}[THEOREM]{Corrolary}

\newtheorem{Subclaim}[THEOREM]{Subclaim}

\newtheorem{Corollary}[THEOREM]{Corollary}

\newtheorem{Example}[THEOREM]{Example}

\newtheorem{Proposition}[THEOREM]{Proposition}

\newtheorem{Discussion}[THEOREM]{Discussion}

\newenvironment{Proof of the Subfact}
{\noindent{\bf Proof of the Subfact.}}{\par\bigskip}

\newenvironment{Proof of the Theorem}
{\noindent{\bf Proof of the Theorem.}}{\par\bigskip}

\newenvironment{Proof of the Conclusion}
{\noindent{\bf Proof of the Conclusion.}}{\par\bigskip}

\newenvironment{Proof of the Observation}
{\noindent{\bf Proof of the Observation.}}{\par\bigskip}

\newenvironment{Proof of the Fact}
{\noindent{\bf Proof of the Fact.}}{\par\bigskip}

\newenvironment{Proof of the Lemma}
{\noindent{\bf Proof of the Lemma.}}{\par\bigskip}

\newenvironment{Proof of the Claim}
{\noindent{\bf Proof of the Claim.}}{\par\bigskip}

\newenvironment{Proof of the Subclaim}
{\noindent{\bf Proof of the Subclaim.}}{\par\medskip}

\newenvironment{Proof of the Main Claim}
{\noindent{\bf Proof of the Main Claim.}}{\par\bigskip}

\catcode`\@=11
\def\@begintheorem#1#2{\rm \trivlist \item[\hskip \labelsep{\bf #1\ #2.}]}
\def\@opargbegintheorem#1#2#3{\rm \trivlist
      \item[\hskip \labelsep{\bf #1\ #2\ (#3).}]}
\catcode`\@=12

\newcommand{\elementary}{\prec}

\newcommand{\Bbf}{\Bbb}

\newcommand{\tensor}{\otimes}

\newcommand{\into}{\rightarrow}

\newcommand{\rest}{\upharpoonright}  

\newcommand{\deq}{\buildrel{\rm def}\over =}

\newcommand{\HH}{{\cal H}}
\newcommand{\II}{{\cal I}}

\newcommand{\PP}{{\cal P}}

\newcount\skewfactor
\def\mathunderaccent#1#2 {\let\theaccent#1\skewfactor#2
\mathpalette\putaccentunder}
\def\putaccentunder#1#2{\oalign{$#1#2$\crcr\hidewidth
\vbox to.2ex{\hbox{$#1\skew\skewfactor\theaccent{}$}\vss}\hidewidth}}
\def\name{\mathunderaccent\tilde-3 }

\begin{document}

\title{Weak reflection at the successor of singular}

\author{Mirna D\v zamonja\\
School of Mathematics\\
University of East Anglia\\
Norwich, NR47TJ, UK\\
\scriptsize{M.Dzamonja@uea.ac.uk}\\
\scriptsize{http://www.mth.uea.ac.uk/people/md.html}\\
Saharon Shelah\\
Mathematics Department\\
Hebrew University of Jerusalem\\
91904 Givat Ram, Israel\\
and\\
Rutgers University\\
New Brunswick, New Jersey\\
USA\\
\scriptsize{shelah@sunset.huji.ac.il}\\
\scriptsize{http://www.math.rutgers.edu/$\thicksim$ shelarch}}

\maketitle
\begin{abstract} {The notion of stationary reflection is one of the
most important notions of combinatorial set theory. We investigate
weak reflection, which is, as the name suggests, a weak version of
stationary reflection. This sort of reflection was introduced in \cite{DjSh
545}, where it was shown that weak reflection has applications
to various guessing principles,
in the sense that if there is no weak reflection, than a guessing principle
holds, and an application dealing with the saturation of normal filters.
Further investigations of weak reflection were carried in  \cite{CuDjSh 571}
and \cite{CuSh xx}. While various $ZFC$ restrictions on the one hand, and 
independence results on the other, were discovered about the weak reflection,
the question of the relative consistency of the existence of a regular
cardinal $\kappa$ such that the first cardinal weakly reflecting at $\kappa$
is a successor of singular, remained open. 
This paper answers that question by proving that (modulo large cardinal
assumptions close to 2-hugeness) that there indeed can be such a
cardinal $\kappa$. {\footnote{This paper is numbered 691 (10/98) in Saharon Shelah's list of
publications. Both authors thank NSF
for partial support by their grant number NSF-DMS-97-04477, as well
as
the United States-Israel Binational Science Foundation for a partial
support through a BSF grant.


}}}
\end{abstract}

\baselineskip=16pt
\binoppenalty=10000
\relpenalty=10000
\raggedbottom

\section{Introduction and the statement of the results.} Stationary
reflection is a compactness phenomenon in the context of stationary sets.
To motivate its investigation,
let us consider first the situation of a regular
uncountable cardinal $\kappa$ and a closed unbounded subset $C$ of $\kappa$.
For every of $\kappa$ many limit points $\alpha$ of $C$, we have that
$C\cap\alpha$ is closed unbounded in $\alpha$. Now let us ask the same
question, but starting with a set $S$ which is stationary, not
necessarily club, in $\kappa$. Is there necessarily $\alpha<\kappa$ such that
$S\cap\alpha$ is stationary in $\alpha$-in the lingo of set theorists,
$S$ reflects at $\alpha$? The answer to this question turns to be very
intricate, and in fact the notion of
stationary reflection is one of the most studied notions of combinatorial set
theory. This is the case not only
because of the historical significance stationary reflection achieved through
by now classical work of R. Jensen \cite{Je} and later work of J.E.
Baumgartner \cite{Baum}, L. Harrington and S. Shelah \cite{HaSh99}, M. Magidor
\cite{Ma} and many later papers, but
also because of the large number of applications it has in set theory and allied areas. In
set theory, stationary reflection is known to have deep connections with
various guessing principles, the simplest one of which is Jensen's $\square$
(\cite{Je}), and the notions from pcf theory, such as good scales (for a long
list of results in this area, as well as an excellent list of references, we
refer the reader to \cite{triple}), and some connections with saturation of
normal filters (\cite{DjSh 545}). In set-theoretic topology, various kinds of
spaces have been constructed from the assumption of the existence of a
non-reflecting stationary set (for references see \cite{handbook}), and in
model theory versions of stationary reflection have been shown to have a connection
with  decidability of monadic second-order logic (\cite{Sh 80}).

We investigate the notion of weak reflection, which, as the name suggests,
is a weakening of the stationary reflection. 
For a regular
cardinal $\kappa$, we say that $\lambda>\kappa$ weakly reflects at $\kappa$
iff for every function $f:\,\lambda\into\kappa$, there is $\delta <\lambda$
of cofinality $\kappa$ (we say $\delta\in
S^\lambda_\kappa$) such that $f\rest e$ is not strictly increasing for any
$e$ a club of $\delta$. 
Its negation is a strong form of non-reflection, called
strong non-reflection. The notions were introduced by D\v zamonja and
Shelah in \cite{DjSh 545} in connection with saturation of normal filters,
as well as
the guessing principle $\clubsuit^\ast_{-\lambda}(\lambda^+)$,
which a relative of another popular guessing
principle, $\clubsuit$.
It is proved in \cite{DjSh 545} that, in the case when
$\lambda=\mu^+$ and $\aleph_0<\kappa=\cf(\mu)<\mu$,
\underline{if}
weak reflection of $\lambda$ at $\kappa$ holds relativized to every
stationary subset of $S^\lambda_\kappa$, \underline{then}
$\clubsuit^\ast_{-\mu}(
S^\lambda_\kappa)$ holds. The exact statement of the principle
is of no consequence to us here, so we omit the definition.
We simply note that this statement is stronger than just $\clubsuit^\ast_{-\mu}
(\lambda)$, which holds just from the given cardinal assumptions.

Weak reflection was further investigated by Cummings, D\v zamonja and Shelah in
\cite{CuDjSh 571}, 
more about which will be mentioned in a moment. Weak reflection
has a very interesting aplication given by Cummings and Shelah in
\cite{CuSh xx}, where they use it as a tool to
build models where stationary reflection holds for some cofinalities but fails
badly for others.

It was proved in \cite{DjSh 545} that if there is $\lambda$ which weakly reflects
at $\kappa$, the first such $\lambda$ is a regular cardinal. It is also not
difficult to see that the first $\lambda$ cannot be weakly compact. On the
other hand, in \cite{CuDjSh 571} Cummings, D\v zamonja and Shelah proved that,
modulo the existence of certain large cardinals, it is consistent to have a
cardinal $\lambda$ which weakly reflects at unboundedly many regular $\kappa$
below it, and strongly non-reflects at unboundedly many others.

A question left open by these investigations, was if it is consistent to have
$\kappa$ for which the first $\lambda$ which weakly reflects at $\kappa$, is a
successor of singular. We answer this question positively, modulo the existence
of a certain large cardinal, whose strength is in the neighborhood of being
2-huge. To state our results more precisely, let us give the exact definition
of weak reflection and the statement of our main theorem.

\begin{Definition}\label{defwr}
Given $\aleph_0<\kappa=\cf(\kappa)$ and $\lambda>\kappa$.
We say that $\lambda$ {\em weakly reflects at} $\kappa$ iff for
every function
$f:\,\lambda\into\kappa$, there is $\delta\in 
S^\lambda_\kappa$ such that $f\rest e$ is not strictly increasing for any
$e$ a club of $\delta$.
\end{Definition}

\begin{Theorem}\label{wref} (1) Let $V$ be a universe in which $GCH$ holds
and $\mu_0$ is a cardinal such that there is an elementary embedding
${\bf j}:\,V\into M$ with the following properties:
\begin{description}
\item{(i)} $\crit({\bf j})=\mu_0$,
\item{(ii)} For some $\kappa^\ast$ a successor of singular and $\chi$, we have
\[ 
\mu_0<\kappa^\ast<\mu_1\deq {\bf j}(\mu_0)<\lambda^\ast\deq {\bf j}(\kappa^\ast)<
\cf(\chi)=\chi<\mu_2\deq{\bf j}(\mu_1),
\]
\item{(iii)} ${}^\chi M\subseteq M$.
\end{description}

\underline{Then} there is a generic extension of $V$ in which cardinals
and cofinalities $\ge \mu_0$ are preserved, and
the first $\lambda$ weakly reflecting at $\kappa^\ast$ is $\lambda^\ast$
(hence, a successor of singular).

{\noindent (2)} In (1), we can replace
the requirement that $\kappa^\ast$ is a successor of singular by
``$\varphi(\kappa^\ast)$ holds"
for any of the following meanings of $\varphi(x)$:
\begin{description}
\item{(a)} $x$ is inaccessible,
\item{(b)} $x$ is strongly inaccessible,
\item{(c)} $x$ is Mahlo,
\item{(d)} $x$ is strongly Mahlo,
\item{(e)} $x$ is $\alpha$-(strongly) inaccessible for $\alpha<x$,
\item{(e)} $x$ is $\alpha$-(strongly) Mahlo for $\alpha<x$,
\end{description}
and have the same conclusion (hence in place of $\lambda^\ast$ is a successor
of
singular, $V^P$ will satisfy $\varphi(\lambda^\ast)$).

{\noindent (3)} With the same assumptions as in (1),
\begin{description}
\item{(i)}
there is a generic extension
of $V$ in which $\kappa^\ast=\aleph_{53}$ and $\lambda^\ast$, a successor of
singular,
is the first cardinal
weakly reflecting at $\kappa^\ast$, 

\item{(ii)} there is a generic extension
of $V$ in which $\kappa^\ast=\aleph_{\omega+1}$ and $\lambda^\ast$, a successor of
singular, is the first
cardinal
weakly reflecting at $\kappa^\ast$, 

\item{(iii)} there is a generic extension
of $V$ in which $\lambda^\ast$ is the first
cardinal
weakly reflecting at $\kappa^\ast$, and $\lambda^\ast=(\kappa^\ast)^{+\gamma}$
for some $\gamma\le\aleph_7$. 

\end{description}
\end{Theorem}

\begin{Remark} Our assumptions follow if $\mu_0$ is 2-huge and $GCH$ holds.
The integers 53 and 7 in the statement of part (3) above, are to a large extent
arbitrary.
\end{Remark}

The proof of (1) uses 
as a building block a forcing notion introduced by Cummings, D\v zamonja and
Shelah in \cite{CuDjSh 571}, which introduces a function witnessing strong
non-reflection of a given cardinal $\lambda$ to a cardinal $\kappa$. An
important feature of this forcing is that it has a reasonable degree of
(strategic) closure, provided that strong non-reflection of $\theta$ to
$\kappa$ already holds for $\theta\in [\kappa,\lambda)$, and hence it can
be iterated. This
forcing is  a rather homogeneous forcing, so the term forcing associated with it has
strong decision properties. The forcing that we actually use is a term forcing
associated with a certain product of the strong non-reflection forcings and
a Laver-like preparation. Using this, we force the strong non-reflection 
of $\theta$ to $\kappa^\ast$ for all $\theta<\lambda^\ast$, and the point is
to prove that in the extension $\lambda^\ast$ weakly reflects on $\kappa^\ast$.
If we are given a condition and a name forced to be a strogly non-reflecting
function, we can use
the large cardinal assumptions to pick a certain model $N$, for which are able
to build a generic condition, whose existence contradicts the choice of the
name. To build the generic condition we use the preparation and the fact that
we are dealing with a term forcing.
Proofs of (2) and (3) are easy modifications of the proof of (1). 

We recall some facts and definitions.

\begin{Notation} (1) Reg stands for the class of regular cardinals.

{\noindent (2)}
If $p,q$ are elements of a forcing notion, then $p\le q$ means
that $q$ is an extension of $p$.

{\noindent (3)} For $p$ a condition in the limit
of an iteration $\langle P_\alpha,\name{Q}_\beta:\,\alpha\le\alpha^\ast,\beta<
\alpha^\ast\rangle$, we let
\[
\Dom(p)\deq\{\beta<\alpha^\ast:\,\neg(p\rest\beta\forces
``p(\beta)=\emptyset_{\name{Q}_\beta}")\}.
\]

{\noindent (4)} The statement that $\lambda$ weakly reflects at $\kappa$
is denoted by $WR(\lambda,\kappa)$.
Its negation  of $WR(\lambda,\kappa)$ is denoted by
$SNR(\lambda,\kappa)$.

\end{Notation}

\begin{Remark}\label{opaska}
It is easily seen that $\lambda$ weakly reflects at $\kappa$
iff $\card{\lambda}$ does, so we can without loss of generality, when
discussing weak reflection of $\lambda$ to $\kappa$ assume that $\lambda$ is a
cardinal. \end{Remark}

\begin{Definition} (1) For a forcing notion and a limit ordinal $\varepsilon$, we
define the game $G(P,
\alpha)$ as follows. The game is played between I and II, and it lasts
$\varepsilon$ steps, unless a
player is forced to stop before that time. 
For $\zeta<\varepsilon$, we denote the $\zeta$-th move
of I by $p_\zeta$, and that of II by $q_\zeta$. The requirements are that
I commences by $\emptyset_P$ and that for all $\zeta$ we have
$p_\zeta\le q_\zeta$, while for $\xi<\zeta$ we have $q_\xi\le p_\zeta$. 

I wins a play $\Gamma$ of $G(P,
\varepsilon)$ iff $\Gamma$ lasts $\varepsilon$ steps.

{\noindent (2)} For $P$ and $\varepsilon$ as above, we say that
$P$ is $\varepsilon$-{\em strategically closed} iff I has a winning
strategy in $G(P,\varepsilon)$. We say that $P$ is
$(<\varepsilon)$-{\em strategically closed} iff it is $\zeta$-strategically
closed for all $\zeta<\varepsilon$.
\end{Definition}

\begin{Definition} A set $A$ of ordinals is an {\em Easton set} iff
\[
\sigma\in {\rm Reg}\cap (\sup(A)+1)\implies\sup(A\cap
\sigma)<\sigma.
\]
\end{Definition}

\begin{Definition} We shall call a forcing notion $P$ {\em mildly homogeneous}
iff for every
formula $\varphi(x_0,\ldots,x_{n-1})$ of the forcing language of $P$ and
$a_0,\ldots,a_{n-1}$ (canonical names of) objects in $V$, we have
$\emptyset_P\decides ``\varphi(a_0,\ldots,a_{n-1})".$
\end{Definition}

\section{Proofs.} We give the proof of Theorem \ref{wref}. The main part of the
proof is to deal with part (1) of the Theorem, and at the very end of the
section we indicate the changes needed for the other parts of the theorem.

\begin{Definition} Given $\aleph_0<\kappa=\cf(\kappa)<\sigma$.

${\Bbf P}(\kappa,\sigma)$ is the forcing notion whose elements are functions
$p$ with $\dom(p)$
an ordinal $<\sigma$, the range $\rge(p) \subseteq\kappa$, and
the
property
\[
\beta\in S^\sigma_\kappa\implies (\exists c \hbox{ a club of }\beta)\,
[p\rest c\mbox{ is strictly increasing}],
\]
while ${\Bbf P}(\kappa,\sigma)$ is ordered by extension.
\end{Definition}

\begin{Fact}[Cummings, D\v zamonja and Shelah]
\label{one} \cite{CuDjSh 571} Let $\kappa$
and $\sigma$ be such that
${\Bbf P}(\kappa,\sigma)$ is defined, \underline{then}
\begin{description}
\item{(1)} $\card{{\Bbf P}(\kappa,\sigma)}\le\card{{}^{<\sigma}\kappa}=
\kappa^{<\sigma}$.
\item{(2)} Suppose that for all $\theta\in [\kappa,\sigma)$ we have
$SNR(\theta,\kappa)$.
\underline{Then} ${\Bbf P}(\kappa,\sigma)$ is $(<\sigma)$-strategically closed.
\end{description}
\end{Fact}

\begin{Definition}\label{iter}
Given $\aleph_0<\cf(\kappa)=\kappa<\lambda$.

$Q_{(\kappa,\lambda)}$ is the result of the reverse Easton support iteration of
${\Bbf P}(\kappa,\sigma)$ for $\sigma=\cf(\sigma)\in (\kappa,\lambda)$.
More precisely, let
\[
\bar{Q}=\langle Q_\alpha,\name{R}_\beta:\,\alpha\le\lambda, \beta<\lambda\rangle,
\]
where 
\begin{description}
\item{(1)} $Q_\alpha\forces``\name{R}_\alpha=\{\emptyset\}"$ unless
$\alpha\in {\rm Reg}\cap(\kappa,\lambda)$, in which case
\[
Q_\alpha\forces``\name{R}_\alpha=\name{{\Bbf P}}(\kappa,\alpha)".
\]
\item{(2)} For $\alpha\le\lambda$

$p\in Q_\alpha$ iff for all $\gamma<\alpha$ we have $\forces_{Q_\gamma}``
p(\gamma)\in \name{R}_\gamma"$ \underline{and}
\begin{description}
\item{(i)}
if $\alpha$ is inaccessible, then
$\card{\Dom(p)}<\alpha$.

\item{(ii)} If $\alpha$ is a limit but not inaccessible, then
\[
Q_\alpha\deq\{p:\,(\forall\beta<\alpha)\,[p\rest\beta\in Q_\beta]\}.
\]
\end{description}

\item{(3)} $p\le q$ iff for all $\beta<\lambda$ we have
$q\rest\beta\forces``q(\beta)\ge p(\beta)"$.
\end{description}

\end{Definition}

\begin{Fact}[Cummings, D\v zamonja, Shelah]\label{closure} \cite{CuDjSh 571}
Let
$\bar{Q}$, $\kappa$ and $\lambda$ be as in Definition \ref{iter}.
For all $\alpha\le\lambda$:
\begin{description}
\item{(1)} Whenever $\alpha$ is regular, $\card{Q_\alpha}\le\alpha^{<\alpha}$,
\item{(2)} $\forces_{Q_\alpha}``\card{
\name{R}_\alpha}\le\card{\alpha}^{<\kappa}"$,
\item{(3)} $Q_{\alpha+1}$ has $(\card{\alpha}^{\card{<\alpha}})^+$-cc. In
addition,
if $\alpha$ is strongly Mahlo, then $Q_\alpha$ has $\alpha$-cc.
\item{(4)} $\forces_{Q_\alpha}``\name{R}_\alpha$ is $(<\alpha)$-strategically
closed".
\item{(5)} For all $\beta<\alpha$, we have that $\name{Q}_\alpha/ Q_\beta$
is $(<\beta)$-strategically closed.
\item{(6)} $Q_\alpha$
preserves all cardinals and cofinalities $\ge
(\card{\alpha}^{<\card{\alpha}})^+$, and
all strongly inaccessible cardinals and cofinalities $\le\card{\alpha}$.
\item{(7)} $\forces_{Q_\alpha}``SNR(\kappa,\beta)"$ for all $\beta<\alpha$.
\end{description}
\end{Fact}

\begin{Notation}
\begin{description}
\item{(1)}
For a forcing notion $Q$ of the form $Q=P_1\ast\name{P}_2$,
we denote by $Q^\tensor$ the term forcing associated with $Q$, defined by
\[
Q^\tensor\deq\{(\emptyset_{P_1},\name{q}):\,\name{q}\mbox{ is a canonical }
P_1\mbox{-name for a condition in }P_2\},
\]
(in particular $Q^\tensor\subseteq Q$), with the order inherited from $Q$.
\item{(2)}
For a triple $(R,\kappa, \sigma)$ with
$\aleph_0<\cf(\kappa)=\kappa<\sigma$,
and $R$ a forcing notion preserving $\kappa=\cf(\kappa)>\aleph_0$,
with $\emptyset_R$ the minimal element of $R$, we define
$Q_{(R,\kappa,\sigma)}^\tensor$
to be $[R\ast\name{Q}_{(\kappa,\sigma)}]^\tensor$.
\end{description}
\end{Notation}

\begin{Observation}\label{closed} $Q_{(R,\kappa,\sigma)}^\tensor$, when defined,
is $(<\kappa^+)$-strategically closed.
\end{Observation}


\begin{Claim}\label{homo}
\begin{description}
\item{(1)} ${\Bbf P}(\sigma,\lambda)$ is mildly homogeneous, for
$\aleph_0<\cf(\sigma)=\sigma
<\lambda$,
\item{(2)} If $P$ is mildly homogeneous and 
\[
\forces_P``\name{Q}\mbox{ is mildly homogeneous}",
\]
\underline{then} $P\ast\name{Q}$ is mildly homogeneous.

\item{(3)} $Q_{(\kappa,\lambda)}$ is mildly homogeneous, for
$\aleph_0<\cf(\kappa)=\kappa<
\lambda$.
\end{description}
\end{Claim}

\begin{Proof of the Claim} (1) Suppose not, and let $p,q\in P\deq{\Bbf
P}(\sigma,\lambda)$
force contradictory statements about $\varphi(a_0,\ldots,a_{n-1})$. Let
$\alpha=\Dom(q)$ and
consider the function $F:\,P\into P$ such that $F(f)=g$ iff $q\subseteq g$ and
for $i\in
\dom(f)$ we have $g(\alpha+i)=f(i)$.

This function is an isomorphism between $P$ and $P/q\deq\{g\in P:
\,g\supseteq q\}$, and it induces an isomorphism between $P$-names
and $P/q$-names. However, in $P/q$ we have that $q$ and $F(p)
\ge q$ force contradictory statements about $\varphi(a_0,\ldots,a_{n-1})$.
Contradiction.

(2)-(3) Similar proofs.
$\eop_{\ref{homo}}$
\end{Proof of the Claim}

\begin{Remark} We remark that ${\Bbf P}(\sigma,\lambda)$ in fact has stronger
homogeneity properties, a fact which will not be used here.
\end{Remark}

{\em Proof of the Theorem continued.}

Let $V$, ${\bf j}$ and the cardinals mentioned in the assumptions
be fixed.
Note that by elementarity, $\lambda^\ast$ is the successor of a singular
cardinal.

\begin{Notation} In the situation when notation
$Q^\tensor_{(R,\kappa,\kappa^\ast)}$
makes sense, we abbreviate it as $Q^\tensor_{(R,\kappa)}$.
\end{Notation}

\begin{Definition} We define ${\Bbf P}^-$ to be the forcing
whose elements are functions $h$, with $\dom(h)$ an Easton subset of
$\mu_0$, with the property 
\[
\alpha<\beta\in \dom(h)
\implies h(\alpha)\in \HH(\beta),
\]
ordered by extension.
\end{Definition}

\begin{Claim}\label{small} Forcing with ${\Bbf P}^-$ preserves cardinals and
cofinalities $\ge\mu_0$, and $GCH$ above $\mu_0$,
while any inaccessible $\sigma\le\mu_0$ which is a limit of inaccessibles,
remains
regular and $2^\sigma=\sigma^+$ holds in the extension by ${\Bbf P}^-$.
\end{Claim}

\begin{Proof of the Claim} First notice that $\card{{\Bbf P}^-}
=\mu_0$, so ${\Bbf P}^-$ has $\mu^+_0$-cc and preserves cardinals and
cofinalities 
$\ge\mu_0^+$,
as well as $GCH$ above $\mu_0$.

Now suppose that $\sigma\le\mu_0$ is a limit of inaccessibles,
but its cofinality is changed by ${\Bbf P}^-$ to be $\le\theta$ for some
$\theta<\sigma$. Let $p^\ast\in {\Bbf P}^-$ force this.
Without loss of generality, $\theta$ is (strongly) inaccessible and $\theta
\in\dom(p^\ast)$.

Let 
\[
P_{<\theta}\deq\{q\rest\theta:\,q\in {\Bbf P}^-\,\,\&\,\,q\ge p^\ast\}
\]
and
\[
P_{\ge\theta}\deq\{q\rest [\theta,\mu_0):\,q\in {\Bbf P}^-\,\,\&\,\,q\ge
p^\ast\}, 
\]
both ordered
by extension. The mapping $q\maps(q\rest[\theta,\mu_0), q\rest\theta)$
shows that ${\Bbf P}^-/p\deq\{q\in {\Bbf P}^-:\,q\ge p^\ast\}$
is isomorphic to $P_{\ge\theta}\times P_{<\theta}$.
We have that $P_{\ge\theta}$ is $(<\theta^+)$-closed, so $P_{<\theta}$
adds a cofinal function from $\theta$ to $\sigma$. However,
$\card{P_{<\theta}}\le\theta$ (as $\theta$ is strongly inaccessible), and
so it preserves cardinals and cofinalities $\ge\theta^+$, a contradiction.

We can similarly decompose ${\Bbf P}^-$ into $P_{\ge\sigma}\times
P_{<\sigma}$ to observe that
\[
\forces_{{\Bbf P}^-}``2^\sigma=\sigma^+".
\]
$\eop_{\ref{small}}$
\end{Proof of the Claim}

\begin{Definition}
{\noindent (1)} For $\mu<\mu_0$ let $\name{{\Bbf R}}_\mu$ and
$\name{\kappa}_\mu$be the following
${\Bbf P}^-$-names:
for a condition $p\in {\Bbf P}^-$, \underline{if} $\mu\in\dom(p)$ and
\[
p(\mu)=(\kappa,\name{R})\mbox{  with }\mu<\cf(\kappa)=\kappa<\mu_0,
\]
and
$R\in \HH(\kappa^+)$ is a
forcing notion which preserves
cardinals and cofinalities $\ge\mu$, \underline{then} $p$ forces
$\name{\kappa}_\mu$ to be $\kappa$ and
$\name{{\Bbf R}}_\mu$
to be $\name{Q}^\tensor_{(R,{\name{\kappa}}_\mu)}$. We say that
$R_\mu=R$.
\underline{If},
$\mu\in\Dom(p)$ but $p(\mu)$ does not satisfy the conditions above,
\underline{then} $p$ forces $\name{{\Bbf R}}_\mu$ to be the trivial forcing,
which will for notational purposes be thought of as $\{\emptyset,\emptyset)\}$.
For the same reason, in these circumstances
we think of $R_\mu=\{\emptyset\}$.

{\scriptsize{Note: each $\name{\Bbf R}_\mu$ is (over a dense subset of
${\Bbf P}^-$) a ${\Bbf P}^-$-name of a forcing notion from $V$,
$\name{\kappa}_\mu$ is a ${\Bbf P}^-$-name of an ordinal $<\mu_0$,
and $\prod_{\mu<\mu_0}
\name{\Bbf R}_\mu$ is a ${\Bbf P}$-name of a product of forcing.
But ${\Bbf R}$ below is forced not to be from $V$.}}

{\noindent (2)} For a ${\Bbf P}^-$-name
$\name{f}\in \prod_{\mu<\mu_0}\name{{\Bbf R}}_\mu$ and
$\alpha\le\mu_0$,
let
\[
\name{A}_{\name{f},\alpha}\deq\{\mu<\mu_0:\,
\name{f}(\mu)=(\emptyset,\name{q})\mbox{ with }
\neg(\forces_{\name{{\Bbf R}}_\mu}``\alpha\notin\Dom(\name{q})")\}.
\]

{\noindent (23)}
Let $\name{{\Bbf R}}$ be a ${\Bbf P}^-$-name for:
\[
\left\{\name{f}\in \prod_{\mu<\mu_0}\name{{\Bbf R}}_\mu:\,(\forall\alpha\le\mu_0)
\,[\name{A}_{\name{f},\alpha}\mbox{ is an Easton set ]}\right\},
\]
ordered by the order inherited from $\prod_{\mu<\mu_0}\name{{\Bbf R}}_\mu$.

{\scriptsize{hence, ${\name{\Bbf R}}$ is a ${\Bbf P}^-$-name of a forcing
notion.}} \end{Definition}

\begin{Notation} 
If we write $(p,\bar{\name{r}})\in {\Bbf P}^-\ast
\name{{\Bbf R}}$, we mean that 
\[
\forces_{{\Bbf P}^-} \bar{\name{r}}=\langle(\emptyset_{R_\mu}, \name{r}
(\mu)):\,\mu<\mu_0\rangle.
\]
\end{Notation}

\begin{Definition} (1) Given $(p,\bar{\name{r}})\in {\Bbf P}^-\ast
\name{{\Bbf R}}$ and $\sigma=\cf(\sigma)<\mu_0$. For
$(q,\bar{\name{s}})\in {\Bbf P}^-\ast
\name{{\Bbf R}}$, we define
\begin{description}
\item{(i)} $(q,\bar{\name{s}})\ge_{{\rm pr},\sigma}(p,\bar{\name{r}})$ iff
\begin{description}
\item{$(\alpha)$} $(q,\bar{\name{s}})\ge(p,\bar{\name{r}})$,
\item{$(\beta)$} $q\rest(\sigma+1)=p\rest(\sigma+1)$,
\item{($\gamma)$} For $\mu<\mu_0$ with $
\neg(q\forces``\name{{\Bbf R}}_\mu$ is trivial"), we
have
\[
(q,\emptyset_{\name{R}_\mu})\forces``
\mbox{if }\name{\kappa}_\mu<\sigma,
\mbox{ then }\name{s}(\mu)\rest (\name{\kappa}_\mu,
\sigma]
=\name{r}(\mu)\rest (\name{\kappa}_\mu, \sigma]".
\]
\end{description}

\item{(ii)} $(q,\bar{\name{s}})\ge_{{\rm apr},\sigma}(p,\bar{\name{r}})$ iff
\begin{description}
\item{$(\alpha)$} $(q,\bar{\name{s}})\ge(p,\bar{\name{r}})$,
\item{($\beta$)} $\dom(q)\cap(\sigma+1,\mu_0)=\dom(p)\cap(\sigma+1,\mu_0)$.
\item{($\gamma)$} For $\mu<\mu_0$ with $
\neg(q\forces``\name{{\Bbf R}}_\mu$ is trivial"), we
have
\[
(q,\emptyset_{\name{R}_\mu})\forces``\name{s}(\mu)\rest (\sigma,\kappa^\ast)
=\name{r}(\mu)\rest (\sigma,\kappa^\ast)".
\]
\end{description}
\end{description}
{\noindent (2)} For $(p,\bar{\name{r}})\in {\Bbf P}^-\ast\name{{\Bbf R}}$ and
$\sigma=\cf(\sigma)\le\mu_0$, we let
\[
Q^-_{(p,\bar{\name{r}}),\sigma}\deq\{(q,\name{\bar{s}}):\,(q,\name{\bar{s}})
\ge_{{\rm apr},\sigma} (p',\bar{\name{r}}')\mbox{ for some }
(p',\bar{\name{r}}')\le_{{\rm pr},\sigma} (p,\bar{\name{r}})\},
\]
ordered as a suborder of ${\Bbf P}^-\ast\name{{\Bbf R}}$.
\end{Definition}

\begin{Claim}\label{intr} Given $(p,\bar{\name{r}})
\le (q,\name{\bar{s}})$ in ${\Bbf P}^-\ast\name{{\Bbf R}}$, and a regular
$\sigma <\mu_0$.

\underline{Then} there is a unique $(t,\name{\bar{z}})$ such that
\[
(p,\bar{\name{r}})\le_{{\rm pr},\sigma}(t,\name{\bar{z}})\le_{{\rm apr},\sigma}
(q,\name{\bar{s}}).
\]
\end{Claim}

\begin{Proof of the Claim} Let $t\deq p\rest (\sigma+1)\cup q\rest (\sigma+1,
\mu_0)$. Hence $t\in {\Bbf P}^-$ and $p\le t\le q$. We define
a ${\Bbf P}^-\ast\name{{\Bbf R}}$-name
$\bar{\name{z}}$ by letting for $\mu<\mu_0$
\[
\name{z}(\mu)\deq\left\{
\begin{array}{ll}
\name{r}(\mu)\rest (\name{\kappa}_\mu,\sigma]\frown\name{q}(\mu)\rest(\sigma,
\kappa^\ast]
&\mbox{ if defined}\\
\name{r}(\mu)
&\mbox{ otherwise.}
\end{array}
\right.
\]
$\eop_{\ref{intr}}$
\end{Proof of the Claim}

\begin{Notation} $(t,\name{\bar{z}})$ as in Claim \ref{intr} is denoted by
intr$\left((p,\bar{\name{r}}),(q,\name{\bar{s}})\right)$.
\end{Notation}

\begin{Claim}\label{strate} For $\sigma=\cf(\sigma)<\mu_0$, the forcing
$({\Bbf P}^-\ast\name{{\Bbf R}}, \le_{{\rm pr},\sigma})$ is
$(<\sigma+1)$-strategically closed.
\end{Claim}

\begin{Proof of the Claim} For every $\mu<\mu_0$ we have
\[
\forces_{{\Bbf P}^-\ast\name{{\Bbf R}}}``\name{{\Bbf R}}_\mu\mbox{ non-trivial}
\implies\name{Q}_{(\name{\kappa}_\mu,\kappa^\ast)}/\name{Q}_{(\name{\kappa}_\mu,
\sigma]}\mbox{ is }(<\sigma+1)\mbox{-strategically closed}."
\]
Hence we can find names $\name{{\rm St}}^\sigma_\mu$ of the corresponding winning
strategies which exemplify the above statement.

Suppose that $\zeta\le\sigma$ and $\langle p_\xi=\langle p^0_\xi,
\name{\bar{p}}^1_\xi
\rangle:\,\xi<\zeta\rangle$, $\langle q_\xi=\langle q^0_\xi, \name{\bar{q}}^1_\xi
\rangle:\,\xi<\zeta\rangle$ are sequences of elements of ${\Bbf P}^-
\ast\name{{\Bbf R}}$
such that
\begin{description}
\item{(1)} For all $\xi<\zeta$ we have $p_\xi\le_{{\rm pr},\sigma}q_\xi$,
\item{(2)} For all $\xi<\zeta$ and $\varepsilon<\xi$ we have
$q_\varepsilon\le_{{\rm pr},\sigma}p_\xi$ and for $\mu<\mu_0$ with
$\neg (p_\xi\forces``\name{{\Bbf R}}_\mu$ is trivial"), we have
$(p_\xi,\emptyset_{\name{R}_\mu}, \name{p}^1_0(\mu)\rest (\name{\kappa}_\mu,
\sigma])\forces_{{\Bbf P}^-\ast\name{{\Bbf R}}_\mu\ast\name{Q}_{(\kappa_\mu,
\sigma]}}$
\[
``
\name{p}^1_\xi(\mu)\rest(\sigma,\kappa^\ast)=\name{{\rm St}}^\sigma_\mu
(\langle\name{p}^1_\varepsilon(\mu)\rest (\sigma,\kappa^\ast):\,\varepsilon
<\xi\rangle,\langle\name{q}^1_\varepsilon(\mu)\rest
(\sigma,\kappa^\ast):\,\varepsilon
<\xi\rangle)"
\]
\end{description}

We define $p_\zeta$ by letting $p^0_\zeta\deq\bigcup\{q^0_\xi:\,\xi<\zeta\}$.
Notice
that $p^0_\zeta\in {\Bbf P}^-$ and $p^0_\zeta\rest
(\sigma+1)=p^0_0\rest(\sigma+1)$.

For $\mu<\mu_0$ with$\neg(p_\xi\forces``\name{{\Bbf R}}_\mu$ is trivial"), we
let
$\name{p}^1_\zeta(\mu)$ be the name given by
\[
\name{p}^1_\zeta(\mu)\rest(\name{\kappa}_\mu,\sigma]\deq
\name{p}^0_\zeta(\mu)\rest(\name{\kappa}_\mu,\sigma]
\]
and
\[
\name{p}^1_\zeta(\mu)\rest(\sigma, \kappa^\ast)\deq
\name{{\rm St}}^\sigma_\mu(\langle\name{p}^1_\xi(\mu)\rest (\sigma,
\kappa^\ast):\,\xi<\zeta\rangle,
\langle\name{q}^1_\xi(\mu)\rest (\sigma,
\kappa^\ast):\,\xi<\zeta\rangle).
\]
$\eop_{\ref{strate}}$
\end{Proof of the Claim}

\begin{Claim}\label{translate}
Suppose $(p,\bar{\name{r}})\forces``\name{\tau}:\,\sigma\into {\rm Ord}"$,
where $\sigma$ is regular $<\mu_0$.

\underline{Then} there is $(q,\bar{\name{s}})\ge_{{\rm pr},\sigma}
(p,\name{\bar{r}})$ and a $Q^-_{(q,\bar{\name{s}}), \sigma}$-name
$\name{\tau}'$ such that 
\[
(q,\bar{\name{s}})\forces``\name{\tau}=\name{\tau}'".
\]
\end{Claim}

\begin{Proof of the Claim} We define a play of $G(({\Bbf P}^-\ast\name{{\Bbf R}},
\le_{{\rm pr},\sigma}),\sigma)$ as follows.

I starts by playing $(p,\bar{\name{r}})\deq p_0$. At the stage $\zeta\le\sigma$,
player
II chooses $q^\ast_\zeta\ge p_\zeta$ such that $q^\ast_\zeta$ forces a value to
$\name{\tau}_\zeta$, and we let $q_\zeta\deq{\rm intr}(p_\zeta, q^\ast_\zeta)$.
At the stage $0<\zeta<\sigma$, we let I play according to the winning strategy
for
$G(({\Bbf P}^-\ast\name{{\Bbf R}},
\le_{{\rm pr},\sigma}),\sigma)$ applied to $(\langle p_\xi:\,\xi<\zeta\rangle,
\langle q_\xi:\,\xi<\zeta\rangle)$. At the end, we let $(q,\bar{\name{s}})
=p_\sigma.$
$\eop_{\ref{translate}}$
\end{Proof of the Claim}

\begin{Claim}\label{cc} If $(p,\bar{\name{r}})\in {\Bbf P}^-\ast\name{{\Bbf R}}$,
and $\sigma=\cf(\sigma)<\mu_0$ is such that $\sigma\in\dom(p)$, \underline{then}
$Q^-_{(p,\name{\bar{r}}),\sigma}$ satisfies $\mu_0$-cc.
\end{Claim}

\begin{Proof of the Claim} Given $\bar{q}=\langle q_i=\langle q^0_i,\name{q}^1_i
\rangle:\,i<\mu_0\rangle$, with $q_i\in Q^-_{(p,\name{\bar{r}}),\sigma}$.
Suppose for contradiction that the range of this sequence is an antichain.

We have that for all $i<\mu_0$
\[
q_i^0\rest (\sigma+1,\mu_0)\subseteq p\rest (\sigma+1,\mu_0).
\]
As $\dom(q^0_i)$ is an Easton set, without loss of generality
we have that all $q_i^0\rest (\sigma+1,\mu_0)$ are the same $q^\ast$. Let $G^-$
be ${\Bbf P}^-$-generic over $V$ with $q^\ast\in G^-$. Hence in $V[G^-]$ the
sequence $\langle q^1_i\deq(\emptyset,\name{\bar{q}}^i):\,i<\mu_0\rangle$
\ to an antichain in $\prod_{\mu<\mu_0}[{\Bbf R}_\mu\ast\name{Q}_{(
\kappa_\mu,\kappa^\ast)}]^\tensor$, and by the choice of the initial sequence,
we have that $\langle(\emptyset, \name{q}^i(\mu)\rest
(\sigma+1)):\,\mu<\mu_0\rangle:\,
i<\mu_0\rangle$ gives rise to an antichain in $\prod_{\mu<\mu_0}[{\Bbf
R}_\mu\ast\name{Q}_{
\kappa_\mu,\sigma}]^\tensor$. For every $i<\mu_0$,
\[
A_i\deq\{\mu<\mu_0:\,\name{q}^i_\mu\rest(\sigma+1)\neq\emptyset\}
\]
has size $\le\sigma$. Hence, without loss of generality, $A_i$'s form a
$\Delta$-system with root $A^\ast$. Hence
\[
\left\langle
\langle\langle\emptyset,\name{q}^i(\mu)\rest(\sigma+1)\rangle:\,\mu\in A^\ast:
\rangle\,i<\mu_0
\right\rangle
\]
gives rise to an antichain, a contradiction.
$\eop_{\ref{cc}}$
\end{Proof of the Claim}

\begin{Claim}\label{preserve}
Forcing with ${\Bbf P}^-\ast\name{{\Bbf R}}$ preserves cardinals
and cofinalities $\ge\mu_0$.
\end{Claim}

\begin{Proof of the Claim} Suppose cofinalities $\ge\mu_0$ are not preserved and
let
$\theta$ be the first
cofinality
$\ge\mu_0$ destroyed. Hence $\theta$ is regular, and for some
$\name{\tau}$, condition $(p, \name{\bar{r}})$
and regular $\sigma<\theta$, we have
$(p, \name{\bar{r}})\forces``\name{\tau}:\,\sigma\into\theta$
is cofinal".

\underline{Case 1}. $\sigma<\mu_0$. By Claim \ref{translate}, there is
$(q,\name{\bar{s}})
\ge (p, \name{\bar{r}})$ and a $Q^-_{(p, \name{\bar{r}}), \sigma}$-name
$\name{\tau}'$ such that
$(q,\name{\bar{s}})\forces``\name{\tau}=\name{\tau}'"$. Hence
$(q,\name{\bar{s}})\forces``\name{\tau}':\,\sigma\into\theta$
is cofinal", contradicting the fact that $Q^-_{(p, \name{\bar{r}}), \sigma}$ has
$\mu_0$-cc.

\underline{Case 2}. $\sigma\ge\mu_0$.

As for every $\mu<\mu_0$ with $\name{{\Bbf R}}_\mu$ non-trivial we have that
\[
{\Bbf P}^-\ast\name{{\Bbf R}}_\mu\ast\name{Q}_{(\name{\kappa}_\mu,\kappa^\ast)}/
\name{Q}_{(\name{\kappa}_\mu,\sigma]}\mbox{ is }(<\sigma+1)
\mbox{-strategically closed},
\]
there is $(q,\name{\bar{s}})\ge (p,\name{\bar{r}})$
and a ${\Bbf P}^-\ast\prod_{\mu<\mu_0}[{\name{\Bbf
R}}_\mu\ast\name{Q}_{(\name{\kappa}_\mu,
\sigma]}]^\tensor$-name $\name{\tau}'$ such that $(q,\name{\bar{s}})\forces
``\name{\tau}':\,\sigma\into\theta$ is cofinal". But this forcing has
$\sigma^+$-cc,
a contradiction.
$\eop_{\ref{preserve}}$
\end{Proof of the Claim}

\begin{Corollary} Forcing with ${\Bbf P}^-\ast\name{{\Bbf R}}
\ast\name{Q}_{(\kappa^\ast,\lambda^\ast)}$
preserves cardinalities and cofinalities $\ge \mu_0$, and forces
$SNR(\theta, \kappa^\ast)$ for $\theta\in (\kappa^\ast,\lambda^\ast)$.
\end{Corollary}

\begin{Claim}\label{homo2} The following is forced by ${\Bbf P}^-$:
\begin{description}
\item{(1)} $\name{{\Bbf R}}$ is mildly homogeneous.
\item{(2)} $\name{{\Bbf R}}\ast\name{Q}_{(\kappa^\ast,\lambda^\ast)}$ is mildly
homogeneous.
\end{description}
\end{Claim}

\begin{Proof of the Claim} (1) First note that each $\name{{\Bbf R}}_\mu$ is
forced to 
be mildly homogeneous, by Claim \ref{homo}(2) and the definition of
$\tensor$ operation.

(2) Follows from (1) and Claim \ref{homo}(2).
$\eop_{\ref{homo2}}$
\end{Proof of the Claim}

\begin{Main Claim}\label{main} After forcing with $
{\Bbf P}\deq{\Bbf P}^-\ast\name{{\Bbf
R}}\ast
\name{Q}_{(\kappa^\ast,\lambda^\ast)}$, we have the weak reflection of
$\lambda^\ast$ at $\kappa^\ast$.
\end{Main Claim}

\begin{Proof of the Main Claim} Suppose otherwise, and let $p^\ast=(p,\name{q},
\name{r})$ force $\name{\tau}$ to be a function exemplifying the strong
non-reflection of $\lambda^\ast$ at $\kappa^\ast$. As $\name{{\Bbf R}}\ast
\name{Q}_{(\kappa^\ast,\lambda^\ast)}$ is forced to be mildly homogeneous
by Claim \ref{homo2},
without
loss of generality $p^\ast=(p,\emptyset,\emptyset)$.

By a standard argument, our large cardinal assumptions imply that
we can find a model $N\elementary\HH(\chi)$ such that
\begin{description}
\item{(i)} $N\cap \mu_0$ is an ordinal $\mu<\mu_0$,
\item{(ii)} $\otp{N\cap\lambda^\ast}=\kappa^\ast$,
\item{(iii)} $\otp{N\cap\mu_1}=\mu_0$,
\item{(iv)} ${}^{\omega}N\subseteq N$ (even ${}^{\mu>}N\subseteq N$,
although we do not use this),
\item{(v)} $(N,\in)$ is isomorphic to $\HH(\chi')$ for some regular
$\chi'<\chi$ \item{(vi)} $\card{N\cap\kappa^\ast}$ is a regular cardinal.
\item{(vii)} $\kappa^\ast,\mu_0,\mu_1,\mu_2,\lambda^\ast, 
{\Bbf P}, p^\ast, \name{\tau}\in N$.
\end{description}

[Why? Consider ${\bf j}``(\HH(\chi))$ and use elementarity. Note that by
${}^\chi M\subseteq M$ and $\chi^{<\chi}=\chi$ we have that ${\bf j}``(\HH(\chi))
\in M$.]

First we consider some consequences of our choice of $N$. Let $\kappa\deq
\card{N\cap\kappa^\ast}$ and $\delta\deq\sup(N\cap\lambda^\ast)$.

Our assumptions on $N$ imply that $\otp{N\cap\kappa^\ast}
<\mu_0$, hence $\kappa<\mu_0$. As for $\delta$, we have
$\delta\in S^{\lambda^\ast}_{\kappa^\ast}$. Now notice that $N\cap\delta$
is stationary in $\delta$, and it remains so after forcing with ${\Bbf P}$.

[Why? The set
$S\deq S^\delta_{\aleph_0}\cap N$ is a stationary subset of $\delta$, as
$E$ defined as the closure of $N\cap\delta$ is a club of $\delta$,
and $[\alpha\in N\,\,\&\,\,\cf(\alpha)=\aleph_0]\implies\alpha\in S$
(this is true even 
with``$\cf(\alpha)<\mu$" in place of ``$\cf(\alpha)=\aleph_0$").
But  ${\Bbf P}$ is an
$\omega_1$-closed forcing notion, hence $S$ remains stationary after forcing with
${\Bbf P}$.]

As
$p^\ast\in N$, we
have that $\dom(p)\subseteq \mu$. Hence 
\[
p^+\deq p\cup
\{\langle \mu, (\kappa, ({\Bbf P}^-\ast\name{{\Bbf R}})^N)
\rangle\}
\]
is a well defined condition in ${\Bbf P}^-$, and it extends $p$. In fact, $p^+$
is a ${\Bbf P}^-$-generic condition over $N$, and it forces that
$\name{{\Bbf R}}_\mu$ is $({\Bbf P}^-\ast\name{\Bbf R})^N$, hence that
$[\name{{\Bbf R}}_\mu\ast\name{Q}_{(\kappa,\kappa^\ast)}]^\tensor$ is
$[({\Bbf P}^-\ast\name{\Bbf R})^N\ast\name{Q}_{(\kappa,\kappa^\ast)}]^\tensor$,
which is $([({\Bbf P}^-\ast\name{{\Bbf
R}})\ast
\name{Q}_{(\kappa^\ast,\lambda^\ast)}]^\tensor)^N$. Let $\name{H}$ be
$([({\Bbf P}^-\ast\name{{\Bbf
R}})\ast
\name{Q}_{(\kappa^\ast,\lambda^\ast)}]^\tensor)^N$-generic with $p\in \name{H}$.
The
inverse of the Mostowski collapse $F$ of $N$ maps $\name{H}$ into a subset
$\name{H}^\ast$ of
$[({\Bbf P}^-\ast\name{{\Bbf
R}})\ast
\name{Q}_{(\kappa^\ast,\lambda^\ast)}]^\tensor$. Notice then that
\[
p^+\forces``\name{H}\mbox{ is }\name{{\Bbf R}}_\mu\mbox{-generic}".
\]
We wish to define $q$ as follows: $q\deq (p^+,\emptyset_{\Bbf R},\name{r})$,
where
$\name{r}$ is a ${\Bbf P}^-\ast\name{{\Bbf R}}$-name over $(p^+,
\emptyset_{\name{\Bbf R}})$ of a condition in
$\name{Q}_{(\kappa^\ast,\lambda^\ast)}$
defined by letting
\[
\Dom(\name{r})\deq\bigcup\{\Dom(\name{h}):\,(p^+,\emptyset_{\name{\Bbf R}})
\forces``F^{-1}(\name{h})\in\name{H}"\},
\]
and for $\theta$ with $(p^+,\emptyset_{\name{\Bbf R}})\forces``\theta
\in\Dom(\name{r})"$, we let
\[
\name{r}(\theta)\deq\bigcup\{\name{h}(\theta):\,(p^+,\emptyset_{\name{\Bbf R}})
\forces``\theta\in\Dom(\name{h})\,\,\&\,\,\name{h}\in\name{H}^\ast"\}.
\]
We now claim that $q$ is a condition in ${\Bbf P}$ and $q\ge p^\ast$.
Let us check  the relevant items:
\begin{description}
\item{(a)} If $(p^+,\emptyset)\forces``\theta\mbox{ strongly inaccessible }
\in(\kappa^\ast,\lambda^\ast)$", then 
\[
(p^+,\emptyset)\forces``
\card{\Dom(\name{r})\cap\theta}<\theta".
\]
\end{description}

[Why? We have that for some $\theta'\in (\kappa,\kappa^\ast)$,
\[
(p^+,\emptyset)\forces``\Dom(\name{r})\cap \theta\subseteq\bigcup\{\Dom
(F(\name{f}):\,\name{f}\in \name{Q}^N_{(\kappa,\theta')}\}".
\]
But then, $(p^+,\emptyset)\forces``\card{\name{Q}^N_{(\kappa,\theta')}}
\le\card{2^{\theta'}\cap N}<\kappa^\ast"$.]

\begin{description}
\item{(b)} If  $(p^+,\emptyset)\forces``\theta\in\Dom(\name{r})"$, then
 $(p^+,\emptyset)\forces``\name{r}(\theta)$ is a function whose domain is
 an ordinal $<\theta$ and range a subset of $\kappa^\ast$".
 
\end{description}

[Why? As $(p^+,\emptyset)\forces``\name{H}^\ast$ is directed", we have that
 $(p^+,\emptyset)\forces``\name{r}(\theta)$ is a function". If
  \[
  (p^+,\emptyset)\forces``\theta\in\Dom(\name{h})\,\,\&\,\,F^{-1}(\name{h})
  \in\name{H}",
  \]
 then $(p^+,\emptyset) $ forces
 \[
 ``(\forall\sigma\in\Dom(F^{-1}(\name{h})))\,
 [F^{-1}(\name{h})(\sigma)\mbox{ is a function with domain }\in 
 \sigma]",
 \]
 so by elementarity
 \[
 (p^+,\emptyset)\forces``\dom(\name{h}(\theta))\mbox{ is an element of
 }\theta".]
 \]
 
 \begin{description}
 \item{(c)} $(p^+,\emptyset)\forces``\name{r}\in
\name{Q}_{(\kappa^\ast,\lambda^\ast)}$.
 \end{description}
 
 [Why? First of all, we know that $(p^+,\emptyset)\forces``\dom(\name{r})
 \subseteq\lambda^\ast$ with no last element". Now,
 for relevant $\theta$, $\dom(\name{r}(\theta))$ is
 the union of a subset of $\HH(\theta^{++})\cap N$ which has cardinality
 $\le\card{\theta^{++}\cap N}<\kappa^\ast$. Hence the union has
 cofinality $<\kappa^\ast$ (as having cofinality $\ge\kappa^\ast$ is
 preserved by the forcing), hence $\name{r}(\theta)$ is forced to be
 in $\name{{\Bbf P}}(\theta,\kappa^\ast)$, and hence $\name{r}$ is forced to
 be an element of $\name{Q}_{(\kappa^\ast,\lambda^\ast)}$.]
 
\begin{description}
\item{(d)}
$ (p^+,\emptyset)$ forces that for all $\varepsilon\in
 S^\theta_{\kappa^\ast}$, there is a club $e$ of $\varepsilon$ on which
$\name{r}(\theta)$
 is strictly increasing, for $\theta\in\dom(\name{r})$.
 \end{description}
 
 [Why? Because this is forced about each $\name{h}(\theta)$ for $\name{h}
 \in\name{H}^\ast$ and there are $<\kappa^\ast$ such $\name{h}$.]
 
 But now we shall see that $q$ forces $\name{\tau}$ to be constant on a
stationary subset
 of $\delta$, a contradiction, as $\delta\in S^{\lambda^\ast}_{\kappa^\ast}$, and
remains
 there after forcing with ${\Bbf P}$. We need to consider what $q$ forces about
$\name{\tau}
 (\alpha)$ for $\alpha\in N$. Such $\name{\tau}(\alpha)$ is a ${\Bbf P}$-name of
an ordinal
 $<\kappa^\ast$. Let
 \[
 \II_\alpha\deq\{(\emptyset,\emptyset,\name{t})\in {\Bbf
P}:\,(\emptyset,\emptyset,\name{t})
 \mbox{ forces }\name{\tau}(\alpha)\mbox{ to be equal to a }
 {\Bbf P}^-\ast\name{{\Bbf R}}\mbox{-name}\}.
 \]
 Hence $\II_\alpha\in N$. As ${\Bbf P}^-\ast\name{{\Bbf R}}$ forces that
$\name{Q}_{(\kappa^\ast,
 \lambda^\ast)}$ is $(\kappa^\ast+1)$-strategically closed, 
 we have that $\II_\alpha$ is dense in $[({\Bbf P}^-\ast
 \name{{\Bbf R}})\ast\name{Q}_{(\kappa^\ast,
 \lambda^\ast)}]^\tensor$. By the definition of $\name{H}^\ast$, 
there is
$(\emptyset,\emptyset,\name{h})\in \II_\alpha\cap N$ such that 
$(\emptyset,\emptyset,\name{h})\le q $. Let $\name{\tau}'$ exemplify this,
so $\name{\tau}'\in N$.

Hence $q$ forces $\name{\tau}(\alpha)$ to be in the set  of all $\name{\tau}'
\in N$,
where $\name{\tau}'$ is a ${\Bbf P}^-\ast\name{{\Bbf R}}$-name of an ordinal
$<\kappa^\ast$.
The cardinality of this set $\le\card{\PP(\kappa^\ast)\cap N}$, which is
$<\mu_0$.
Since $\alpha\in N$ was arbitrary, $q$ forces the range of $\name{\tau}\rest
(N\cap\delta)$ to be a set of size $<\mu_0$, hence $\name{\tau}$ will be constant
on a stationary subset of $N\cap\delta$ (as  $N\cap\delta$ is stationary).
$\eop_{\ref{main}}$
\end{Proof of the Main Claim}

{\em Proof of the Theorem continued.}

(2) Same proof.

(3) Follow the forcing from (1) by a Levy collapse.
We are making use of the following 

\begin{Claim}\label{collapse} Suppose $\lambda$ weakly reflects at $\kappa$
and $P$ is a $\kappa$-cc forcing.

\underline{Then} $\lambda$ weakly reflects at $\kappa$ in $V^P$.
\end{Claim}

\begin{Proof of the Claim} Suppose that 
\[
p\forces_P``\name{f}:\,\lambda\into\kappa".
\]
We define $f':\,\lambda\into\kappa$ by letting $f'(\alpha)\deq
\sup\{\gamma<\kappa:\,\neg (p\forces``\name{f}(\alpha)\neq\gamma")\}$. As
$P$ is $\kappa$-cc, the range of $f'$ is indeed contained in $\kappa$.
Let $\delta\in S^\lambda_\kappa$ be such
that $f'\rest S$ is constant on a stationary set $S\subseteq\delta$ (the
existence of such $\delta$ follows as $WR(\lambda,\kappa)$ holds).
Hence $p\forces``\name{f}\rest S$ is bounded", so $\name{f}$
does not witness $SNR(\lambda,\kappa)$ in $V^P$, as $S$ remains stationary in 
$V^P$.
$\eop{\ref{collapse}}$
\end{Proof of the Claim}

So, for example to get $\kappa^\ast=\aleph_n$ for $n\ge 1$, we could in
$V^{\Bbf P}$ from (1) first make $GCH$ hold below $\mu_0$ by collapsing various cardinals
below $\mu_0$, and then collapse $\kappa^\ast$ to be $\aleph_n$.

If we start in $V$ by having $\kappa^\ast=\mu_0^{\omega+1}$, and in $V^{\Bbf
P}$ collapse $\mu_0$ to $\aleph_1$, then we get $\kappa^\ast=\aleph_{\omega+1}$.

For the last statement, start
by $\kappa^\ast=\mu_0^{+\delta}$ for some $\delta\le \aleph_7$,
so $\lambda^\ast=\mu_1^{+\delta}$.
By a minor change in the definition of ${\Bbf P}^-$,
we can make ${\Bbf P}$ not add any $\aleph_7$-sequences from $V$. Now in
$V^{\Bbf P}$, collapse first $\mu_0$ to $\aleph_7$, and then collapse
$\mu_1$ to $(\kappa^\ast)^{+3}$. Hence $\lambda^\ast=(\kappa^\ast)^{+3+\delta}$.
$\eop_{\ref{wref}}$

\eject


\begin{thebibliography}{xxxxxxx}

\bibitem[Ba]{Baum} J.E. Baumgartner, {\em A new class of order types}, Annals
of Mathematical Logic 9 (1976), 187-222.

\bibitem[CuFoMa]{triple} J. Cummings, M. Foreman and M. Magidor, {\em
Squares, scales and stationary reflection}, to appear.

\bibitem[CuSh 596]{CuSh xx}
J. Cummings and S. Shelah,
{\em Some independence results on reflection},
Journal of the London Mathematical Society 59 (1999) 37-49. 

\bibitem[CuDjSh 571]{CuDjSh 571} J.Cummings, M. D\v zamonja and S. Shelah,
{\em A consistency result on weak reflection}, Fundamenta Mathematicae
148 (1995), 91-100.

\bibitem[DjSh 545]{DjSh 545} M. D\v zamonja and S. Shelah,
{\em Saturated filters at successors of singulars, weak reflection
and yet another weak club principle}, Annals of Pure And Applied Logic 79
(1996) 289-316.

\bibitem[HaSh 99]{HaSh99} L. Harrington and S. Shelah, {\em Some exact
equiconsistency results in set theory}, Notre Dame Journal of Formal Logic,
1985.

\bibitem[Je]{Je} R.B. Jensen, {\em The fine structure of the
constructible hierarchy}, Annals
of Mathematical Logic 4 (1972), 229-308.

\bibitem[KuVa]{handbook} K. Kunen and J. Vaughan (eds.(, {\em Handbook of
Set-theoretic Topology}, North-Holland 1983.

\bibitem[La]{Laver} R. Laver, {\em Making the supercompactness of
$\kappa$ indestructible under $\kappa$-directed closed forcing},
Israel Journal of Math, 29 (1978), 385-388.

\bibitem[Ma]{Ma} M. Magidor, {\em Reflecting Stationary Sets},
Journal of Symbolic Logic 47 No. 4 (1982) 755-771.

\bibitem[Sh 80]{Sh 80} S. Shelah, {\em A weak generalization of $MA$ to higher
cardinlas}, Israel Journal of Mathematics, 30 (1978), 297-360.

\bibitem[Sh 667]{Sh 667} S. Shelah, {\em Not collapsing cardinals $\le\kappa$
in $(<\kappa)$-support iterations II}, in preparation.

\end{thebibliography}
\end{document}